\documentclass[11pt]{article}
\usepackage[english]{babel}
\usepackage[latin1]{inputenc}
\usepackage{exscale}
\usepackage[centertags]{amsmath}
\usepackage{amsfonts,amstext,amssymb,amsthm}
\usepackage{comment}
\usepackage{newlfont}
\usepackage{enumerate}
\usepackage{graphicx}
\usepackage{these}
\usepackage{color}
\usepackage{hyperref}
\usepackage{url}
\title{Generalized compactness for finite perimeter sets and applications to the isoperimetric problem}
\author{Abraham Mu\~noz Flores\footnote{Partially supported by Capes}, Stefano Nardulli}
\begin{document}
      \maketitle
      \begin{center}
\noindent {\sc abstract}. For a complete noncompact connected Riemannian manifold with bounded geometry, we prove a compactness result for sequences of finite perimeter sets with uniformly bounded volume and perimeter in a larger space obtained by adding limit manifolds at infinity. We extends previous results contained in \cite{NarAsian}, in such a way that the generalized existence theorem, Theorem $1$ of \cite{NarAsian} is actually a generalized compactness theorem. The suitable modifications to the arguments and statements of the results in \cite{NarAsian} are non-trivial. As a consequence we give a multipointed version of Theorem 1.1 of \cite{LangWenger}, and a simple proof of the continuity of the isoperimetric profile function.
\bigskip\bigskip

\noindent{\it Key Words:} Existence of isoperimetric region, isoperimetric profile.
\bigskip

\centerline{\bf AMS subject classification: }
49Q20, 58E99, 53A10, 49Q05.
\end{center}
      \tableofcontents       
      \newpage
\section{Introduction}\label{1}
In this paper we study compactness of the set of finite perimeter sets with bounded volume and bounded perimeter, but with possibly unbounded diameter in a complete noncompact Riemannian manifold, assuming some bounded geometry conditions on the ambient manifold. The difficulty is that for a sequence of regions some volume may disappear to infinity. Given a sequence of finite perimeter sets inside a manifold with bounded geometry, we show that up to a subsequences the original sequence splits into an at most countable number of pieces which carry a positive fraction of the volume, one of them possibly staying at finite distance and the others concentrating along divergent directions. Moreover, each of these pieces will converge to a finite perimeter set lying in some pointed limit manifold, possibly different from the original. So a limit finite perimeter set exists  in a generalized multipointed convergence. The range of applications of these results is wide. The vague notions invoked in this introductory paragraph will be made clear and rigorous in what follows. 
\subsection{Finite perimeter sets and compactness}
In the remaining part of this paper we always assume that all the Riemannian manifolds $M$ considered are smooth with smooth Riemannian metric $g$. We denote by $V_g$ the canonical Riemannian measure induced on $M$ by $g$, and by $A_g$ the $(n-1)$-Hausdorff measure associated to the canonical Riemannian length space metric $d$ of $M$. When it is already clear from the context, explicit mention of the metric $g$ will be suppressed in what follows.  
\begin{Def} Let $M$ be a Riemannian manifold of dimension $n$, $U\subseteq M$ an open subset, $\mathfrak{X}_c(U)$ the set of smooth vector fields with compact support on $U$. Given $E\subset M$ measurable with respect to the Riemannian measure, the \textbf{perimeter of $E$ in $U$}, $ \mathcal{P}(E, U)\in [0,+\infty]$, is
      \begin{equation}
                 \mathcal{P}(E, U):=sup\left\{\int_{U}\chi_E div_g(X)dV_g: X\in\mathfrak{X}_c(U), ||X||_{\infty}\leq 1\right\},
      \end{equation}  
where $||X||_{\infty}:=\sup\left\{|X_p|_{{g}_p}: p\in M\right\}$ and $|X_p|_{{g}_p}$ is the norm of the vector $X_p$ in the metric $g_p$ on $T_pM$. If $\mathcal{P}(E, U)<+\infty$ for every open set $U$, we call $E$ a \textbf{locally finite perimeter set}. Let us set $\mathcal{P}(E):=\mathcal{P}(E, M)$. Finally, if $\mathcal{P}(E)<+\infty$ we say that \textbf{$E$ is a set of finite perimeter}.    
\end{Def}
\begin{Def} We say that a sequence of finite perimeter sets $E_j$ \textbf{converges in $L^1_{loc}(M)$} to another finite perimeter set $E$, and we denote this by writing $E_j\rightarrow E$ in $L^1_{loc}(M)$, if $\chi_{E_j}\rightarrow\chi_{E}$ in $L^1_{loc}(M)$, i.e., if $V((E_j\Delta E)\cap U)\rightarrow 0\;\forall U\subset\subset M$. Here $\chi_{E}$ means the characteristic function of the set $E$ and the notation $U\subset\subset M$ means that $U\subseteq M$ is open and $\overline{U}$ (the topological closure of $U$) is compact in $M$.
\end{Def}   
\begin{Def}
We say that a sequence of finite perimeter sets $E_j$ \textbf{converge in the sense of finite perimeter sets} to another finite perimeter set $E$, if $E_j\rightarrow E$ in $L^1_{loc}(M)$, and 
\begin{eqnarray*} 
                           \lim_{j\rightarrow+\infty}\mathcal{P}(E_j)=\mathcal{P}(E).
\end{eqnarray*}
\end{Def}
For a more detailed discussion on locally finite perimeter sets and functions of bounded variation on a Riemannian manifold, one can consult \cite{MPPP}.
The following classical compactness theorem is found in the literature in Theorem 1.19 of \cite{Giu}, Theorem 3.23 \cite{AmbrosioFuscoPallara}, just to cite some good texts on the subject of finite perimeter sets or \cite{Morgmt}, \cite{Fed} for the general case of currents in the context of Euclidean spaces. The following Riemannian analogue is an easy consequence of that arguments and the theory of $BV$-functions in the Riemannian context developed in \cite{MPPP} and thus we omit the proof.
\begin{Thm}\label{ClassicalCompactness} Let $(\Omega_i)_{i\in\N}$ be a sequence of finite perimeter sets such that there exists a positive constant $C>0$ and a fixed large geodesic ball $B$, satisfying $V(\Omega_i)+\P(\partial\Omega_i)\leq C$, $\Omega_i\subseteq B$, for every $i\in\N$. Then there exists a finite perimeter set $\Omega\subseteq B$ such that $\Omega_i$ tends to $\Omega$, in $L^1(M)$ topology. 
\end{Thm}
\subsection{Isoperimetric profile and existence of isoperimetric regions}
Typically in the literature, we give the following definition.
\begin{Def}\label{Def:IsPStrong}
The \textbf{isoperimetric profile function} (or briefly, the \textbf{isoperimetric profile}) $I_M:[0,V(M)[\rightarrow [0,+\infty [$, is defined by $$I_M(v):= \inf\{A(\partial \Omega): \Omega\in \tau_M, V(\Omega )=v \}, v\neq 0,$$ and $I_M(0)=0$, where $\tau_M$ denotes the set of relatively compact open subsets of $M$ with smooth boundary.
\end{Def}
However there is a more general context in which to consider this notion that will be better suited to our purposes.  Namely, we replace the set $\tau_M$ with subsets of finite perimeter.
By standard results of the theory of sets of finite perimeter, we have that $A(\partial^* E)=\mathcal{H}^{n-1}(\partial^*E)=\mathcal{P}(E)$ where $\partial^*E$ is the reduced boundary of $E$. In particular, if $E$ has smooth boundary, then $\partial^*E=\partial E$, where $\partial E$ is the topological boundary of $E$. In the sequel we will not distinguish between the topological boundary and reduced boundary when no confusion can arise.
\begin{Def}\label{Def:IsPWeak}
Let $M$ be a Riemannian manifold of dimension $n$ (possibly with infinite volume). We denote by $\tilde{\tau}_M$ the set of  finite perimeter subsets of $M$. The function $\tilde{I}_M:[0,V(M)[\rightarrow [0,+\infty [$  defined by 
     $$\tilde{I}_M(v):= \inf\{\mathcal{P}(\Omega)=A(\partial \Omega): \Omega\in \tilde{\tau}_M, V(\Omega )=v \}$$ 
is called the \textbf{isoperimetric profile function} (or shortly the \textbf{isoperimetric profile}) of the manifold $M$. If there exist a finite perimeter set $\Omega\in\tilde{\tau}_M$ satisfying $V(\Omega)=v$, $\tilde{I}_M(V(\Omega))=A(\partial\Omega)= \mathcal{P}(\Omega)$ such an $\Omega$ will be called an \textbf{isoperimetric region}, and we say that $\tilde{I}_M(v)$ is \textbf{achieved}. 
\end{Def} 
\begin{Thm}[A. Mu\~noz Flores \cite{Abrahamtese}] If $M^n$ is complete then $I_M=\tilde{I}_M$. 
\end{Thm} 
\begin{Dem} 
The proof of this fact is a consequence of the approximation theorem for finite perimeter sets by a sequence of smooth domains, as stated in the context of Riemannian manifolds in Proposition 1.4 of \cite{MPPP}, and of the article \cite{Modica}, that permits to show that an arbitrary finite perimeter set $E$ could be approximated in $L^1_{loc}$-norm by relatively compact domains $\Omega$ with smooth boundary and $V(\Omega)=V(E)$.
\end{Dem}

Compactness implies always existence of isoperimetric regions, but there are examples of manifolds without isoperimetric regions of some or every volumes. For further information about this point the reader could see the introduction of \cite{NarAsian} and the discussion therein. So we cannot have always a compactness theorem if we stay in a noncompact ambient manifolds. If $M$ is compact, classical compactness arguments of geometric measure theory as Theorem \ref{ClassicalCompactness} combined with the direct method of the calculus of variations provide existence of isoperimetric regions in any dimension $n$. For completeness we remind the reader that if $n\leq 7$, then the boundary $\partial\Omega$ of an isoperimetric region is smooth. More generally, the support of the boundary of an isoperimetric region is the disjoint union of a regular part $R$ and a singular part $S$. $R$ is smooth at each of its points and has constant mean curvature, while $S$ has Hausdorff-codimension at least $7$ in $\partial\Omega$. For more details on regularity theory see \cite{Morg1} or \cite{Morgmt} Sect. 8.5, Theorem 12.2.
 
The behavior of sequences of finite perimeter sets with bounded volume and bounded perimeter in noncompact ambient manifolds is caracterized, roughly speaking by a part that stays at finite distance and a part that goes to infinity. 
In this paper, we describe exactly what happens to the divergent part in the case of bounded geometry, which in this context means that both the Ricci curvature and the volume of geodesic balls of a fixed radius are bounded below. For any fixed $v$, we can build an example of a $2$-dimensional Riemannian manifold obtained by altering a Euclidean plane an infinite number (or under stronger assumptions, a finite number) of sequences of caps, diverging in a possibly infinite number of divergent directions. Our main Theorem \ref{Res:MainCompactness} shows that is essentially all that can occur, in bounded geometry and that the limit region could live in a larger space which includes the limits at infinity. In particular if we start with a minimizing sequence of volume $v$ for the isoperimetric profile we obtain in the limit an isoperimetric region which lies in a finite number of limit manifolds, because of the Euclidean type isoperimetric constant for small volumes due to Theorem ...of  Maheux and Saloff Coste \cite{}, recovering Theorem $1$ of \cite{NarAsian}.
\subsection{Main Results} 
The first result of this paper is Theorem \ref{Res:MainCompactness}, which provides a generalized compactness result for isoperimetric regions in a noncompact Riemannian manifold satisfying the condition of bounded geometry. In the general case, limit finite perimeter sets do not exist in the original ambient manifold, but rather in the disjoint union of an at most countable family of pointed limit manifolds $M_{1,\infty},\dots,M_{N,\infty}$, obtained as limit of $N$ sequences of pointed manifolds $(M, p_{ij}, g)_j$, $i\in\{1,\dots, N\}$, here $N$ is allowed to be $\infty$. When we start with a minimizing sequence, then  $N<+\infty$ and we recover Theorem $1$ of \cite{NarAsian}. To prove Theorem  \ref{Res:MainCompactness}, we get a decomposition lemma (Lemma \ref{Lemma: SequenceStructure}) for the thick part of a subsequence of an arbitrary sequence with uniformly bounded volume and area, but that is not contained in a fixed ball. Lemma  \ref{Lemma: SequenceStructure} is interesting for himself.

Now, let us recall the basic definitions from the theory of convergence of manifolds, as exposed in \cite{Pet}. This will help us to state the main result in a precise way. 
\begin{Def} For any $m\in\mathbb{N}$, $\alpha\in [0, 1]$, a sequence of pointed smooth complete Riemannian manifolds is said to \textbf{converge in
the pointed $C^{m,\alpha}$, respectively $C^{m}$ topology to a smooth manifold $M$} (denoted $(M_i, p_i, g_i)\rightarrow (M,p,g)$), if for every $R > 0$ we can find a domain $\Omega_R$ with $B(p,R)\subseteq\Omega_R\subseteq M$, a natural number $\nu_R\in\mathbb{N}$, and $C^{m+1}$ embeddings $F_{i,R}:\Omega_R\rightarrow M_i$, for large $i\geq\nu_R$ such that $B(p_i,R)\subseteq F_{i,R} (\Omega_R)$ and $F_{i,R}^*(g_i)\rightarrow g$ on $\Omega_R$ in the $C^{m,\alpha}$, respectively $C^m$ topology. 
\end{Def}\noindent 
It is easy to see that this type of convergence implies pointed Gromov-Hausdorff convergence. When
all manifolds in question are closed, the maps $F_i$ are global $C^{m+1}$ diffeomorphisms. So for closed manifolds we can speak about unpointed convergence. What follows is the precise definition of \textbf{$C^{m,\alpha}$-norm at scale $r$}, which can be taken as a possible definition of bounded geometry. 
\begin{Def}[\cite{Pet}]\label{Def:BoundedGeometryPetersen}
A subset $A$ of a Riemannian n-manifold $M$  has \textbf{bounded $C^{m,\alpha}$ norm on the scale of $r$}, $||A||_{C^{m,\alpha},r}\leq Q$, if every point p of M lies in an open set $U$ with a chart $\psi$ from the Euclidean $r$-ball into $U$ such that
\begin{enumerate}[(i):]
      \item For all $p\in A$ there exists $U$ such that $B(p,\frac{1}{10}e^{-Q}r)\subseteq U$.
      \item $|D\psi|\leq e^Q$ on $B(0,r)$ and $|D\psi^{-1}|\leq e^Q$ on $U$.
      \item $r^{|j|+\alpha}||D^jg||_{\alpha}\leq Q$ for all multi indices j with $0\leq |j|\leq m$, where $g$ is the matrix of functions of metric coefficients in the $\psi$ coordinates regarded as a matrix on $B(0,r)$. 
\end{enumerate} 
We write that $(M,g,p)\in\mathcal{M}^{m,\alpha}(n, Q, r)$, if $||M||_{C^{m,\alpha},r}\leq Q$. 
\end{Def}

In the sequel, unless otherwise specified, we will make use of the technical assumption on $(M,g,p)\in\mathcal{M}^{m,\alpha}(n, Q, r)$ that $n\geq 2$, $r,Q>0$, $m\geq 1$, $\alpha\in ]0,1]$.  Roughly speaking, $r>0$ is a positive lower bound on the injectivity radius of $M$, i.e., $inj_M>C(n,Q,\alpha,r)>0$.  
\begin{Def}\label{Def:BoundedGeometry}
A complete Riemannian manifold $(M, g)$, is said to have \textbf{bounded geometry} if there exists a constant $k\in\mathbb{R}$, such that $Ric_M\geq k(n-1)$ (i.e., $Ric_M\geq k(n-1)g$ in the sense of quadratic forms) and $V(B_{(M,g)}(p,1))\geq v_0$ for some positive constant $v_0$, where $B_{(M,g)}(p,r)$ is the geodesic ball (or equivalently the metric ball) of $M$ centered at $p$ and of radius $r> 0$.
\end{Def}
\begin{Rem} In general, a lower bound on $Ric_M$ and on the volume of unit balls does not ensure that the pointed limit metric spaces at infinity are still manifolds. 
\end{Rem}
\begin{Rem} We observe here that Definition \ref{Def:BoundedGeometryInfinity} is weaker than Definition \ref{Def:BoundedGeometryPetersen}. In fact, using Theorem $72$ of \cite{Pet}, one can show that if a manifold $M$ has bounded $C^{m,\alpha}$ norm on the scale of $r$ for $\alpha>0$ in the sense of Definition \ref{Def:BoundedGeometryPetersen} then $M$ has $C^{m,\beta}$-locally asymptotic bounded geometry in the sense of Definition \ref{Def:BoundedGeometryInfinity}, for every $0<\beta<\alpha$, while in general the converse is not true. 
\end{Rem} 
This motivates the following definition, that is suitable for most applications to general relativity see for example \cite{NarFloresSchwarzschild}. 
\begin{Def}\label{Def:BoundedGeometryInfinity}
We say that a smooth Riemannian manifold $(M^n, g)$ has $C^{m,\alpha}$-\textbf{locally asymptotic bounded geometry} if it is of bounded geometry and if for every diverging sequence of points $(p_j)$, there exist a subsequence $(p_{{j}_{l}})$ and a pointed smooth manifold $(M_{\infty}, g_{\infty}, p_{\infty})$ with $g_{\infty}$ of class $C^{m,\alpha}$ such that the sequence of pointed manifolds $(M, p_{{j}_{l}}, g)\rightarrow (M_{\infty}, g_{\infty}, p_{\infty})$, in  $C^{m,\alpha}$-topology.  
\end{Def}
For a more detailed discussion about this point the reader could find useful to consult \cite{NarAsian}.  
\begin{Def}[Multipointed Gromov-Hausdorff convergence]
We say that a sequence of multipointed proper metric spaces $(X_i, d_i, p_{i,1}, ..., p_{i,l},...)$ converges to the multipointed metric space $(X_{\infty}, d_{\infty}, p_{\infty,1}, ..., p_{\infty,l},...)$, in the \textbf{multipointed Gromov-Hausdorff topology}, if for every $j$ we have $$(X_i, d_i, p_{i,j})\rightarrow (X_{\infty}, d_{\infty}, p_{\infty,j}),$$ in the pointed Gromov-Hausdorff topology. 
\end{Def}
\begin{Def}[Multipointed $C^0$-convergence]
We say that a sequence of multipointed Riemannian manifolds $(M_i, g_i, p_{i,1}, ..., p_{i,l},...)$ converges to the multipointed Riemannian manifold $(M_{\infty}, g_{\infty}, p_{\infty,1}, ..., p_{\infty,l},...)$, in the \textbf{multipointed $C^0$-topology}, if for every $j$ we have $$(M_i, d_i, p_{i,j})\rightarrow (M_{\infty}, d_{\infty}, p_{\infty,j}),$$ in the pointed $C^0$-pointed topology. 
\end{Def}
\begin{Rem} Multipointed $C^0$-convergence is stronger than Multipointed Gromov-Hausdorff convergence.
\end{Rem}
\begin{Rem} The perimeter is lower semicontinuous with respect to multipointed $C^0$-topology. The volume is continuous with respect to multipointed Gromov-Hausdorff topology. This last assertion about volumes is a deep result due to Tobias Colding \cite{Colding97Ann}.
\end{Rem}
\begin{Res}[Main: Generalized Compactness]\label{Res:MainCompactness} Let $M^n$ be a complete Riemannian manifold with $C^0$-bounded geometry. Let $(\Omega_j)$ be a sequence of finite perimeter sets with $\P(\Omega_j)\leq A$ and $V(\Omega_j)\leq v$. Then there exist a subsequence $(\Omega_{j_k})$ that we rename by $(\Omega_k)$, $k\in S\subseteq\N$, a double sequence of points $p_{ik}\in M^n$, a finite perimeter set $\Omega\subseteq\tilde{M}$ and a sequence $(p_{\infty,i})$ of points of $\tilde{M}$ such that $\left(\Omega_k, \left(p_{ik}\right)\right)\rightarrow (\Omega, (p_{\infty,i}))$ in the multipointed $C^0$-topology.
\end{Res}
\begin{CorRes}[Generalized existence \cite{NarAsian}]\label{Cor: Genexistence}
Let $M$ have $C^0$-locally asymptotically bounded geometry. Given a positive volume $0<v < V(M)$, there are a finite number $N$, of limit manifolds at infinity such that their disjoint union with M contains an isoperimetric region of volume $v$ and perimeter $I_M(v)$. Moreover, the number of limit manifolds is at worst linear in $v$.
Indeed $N\leq\left[\frac{v}{v^*}\right]+1=l(n,k, v_0, v)$, where $v^*$ is as in Lemma 3.2 of \cite{EB}.
\end{CorRes}
\begin{Rem} Observe that if $(M,g,p)\in\mathcal{M}^{m,\alpha}(n, Q, r)$ for every $p\in M$, then $M$ have $C^0$-bounded geometry. So Theorem \ref{Res:MainCompactness} and of course Corollary \ref{Cor: Genexistence} applies to pointed manifolds $M\in\mathcal{M}^{m,\alpha}(n, Q, r)$, hence a posteriori also to manifolds with bounded (from above and below at the same time) sectional curvature and positive injectivity radius. 
\end{Rem}

Now we come back to the main interest of our theory, i.e., to extend arguments valid for compact manifolds to noncompact ones. To this aim let us introduce
the following definition suggested by Theorem \ref{Cor: Genexistence}.
\begin{Def} Let $M$ be a $C^0$-locally asymptotic bounded geometry Riemannian manifolds. We call $D_{\infty}=\bigcup_{i} D_{\infty, i}$ a finite perimeter set in $\tilde{M}$ a \textbf{generalized set of finite perimeter of $\tilde{M}$} and an isoperimetric region of $\tilde{M}$ a \textbf{generalized isoperimetric region}, where $\tilde{M}:=\{(\mathcal{N},q, g_{\mathcal{N}}):\exists (p_j), p_j\in M, p_j\rightarrow +\infty, (M,p_j, g)\rightarrow (\mathcal{N},q, g_{\mathcal{N}})\}$. 
\end{Def}
\begin{Rem} We remark that $D_{\infty}$ is a finite perimeter set in volume $v$ in $\mathring{\bigcup}_i M_{\infty, i}$. \end{Rem}
\begin{Rem} If $D$ is a genuine isoperimetric region contained in $M$, then $D$ is also a generalized isoperimetric region with $N=1$ and $$(M_{\infty,1}, g_{\infty,1})=(M,g).$$ This does not prevent the existence of another generalized isoperimetric region of the same volume having more than one piece at infinity.
\end{Rem}
\subsection{Plan of the article}
\begin{enumerate}
           \item  Section \ref{1} constitutes the introduction of the paper. We state the main results of the paper. 
           \item In Section \ref{2}, we prove Theorem \ref{Res:MainCompactness} as 
           stated in Section \ref{1}.
           \item In the Appendix  we describe some useful links to the existing literature related with Theorem \ref{Res:MainCompactness}. In particular we give a multipointed version of the pointed compactness Theorem 1.1 of \cite{LangWenger} 
\end{enumerate}
\subsection{Acknowledgements}  
The authors would like to aknowledge Christina Sormani, Pierre Pansu, Andrea Mondino, Michael Deutsch, Frank Morgan, Manuel Ritor\'e, Luigi Ambrosio for their useful comments and remarks. The first author wish to thank the CAPES for financial support.
\newpage    
      \section{Proof of Theorem \ref{Res:MainCompactness}}\label{2}
The general strategy used in calculus of variations to understand the structure of solutions of a variational problem in a noncompact ambient manifold is the Concentration-Compactness principle of P.L. Lions.  This principle suggests an investigation of regions in the manifold where volume concentrates.  For the aims of the proof, this point of view it is not strictly necessary.  But we prefer this language because it points the way to further applications of the theory developed here for more general geometric variational problems, PDE's, and the Calculus of Variations.  The concentration compactness Lemma 2.1 of \cite{NarAsian} or Lemma I.1 of \cite{Lions} will be used in our problem, taking measures $\mu_j$ having densities $\chi_{D_j}$,
 where $\chi_{D_j}$ is the characteristic function of $D_j$ which have bounded volumes and bounded perimeter and in particular for an almost minimizing sequence $(D_j)$ defined below. 
\begin{Def}
We say that $(D_j)_j\subseteq\tilde{\tau}_M$ (see Defn. \ref{Def:IsPWeak}) is an \textbf{almost minimizing sequence} in volume $v>0$ if \begin{enumerate}[(i):]
       \item $V(D_j)\rightarrow v$,
       \item $A(\partial D_j)\rightarrow I_M(v)$.
\end{enumerate}
\end{Def} 
\subsection{Proof of Theorem \ref{Res:MainCompactness}.}
 For what follows it will be useful to give the definitions below. 
\begin{Def}
Let $\phi:M\rightarrow N$ a diffeomorphism between two Riemannian manifolds and $\varepsilon>0$. We say that $\phi$ is a \textbf{$(1+\varepsilon)$-isometry} if for every $x,y\in M$, $(1-\varepsilon)d_M(x,y)\leq d_N(\phi(x),\phi(y))\leq (1+\varepsilon)d_M(x,y)$.       
\end{Def}
For the reader's convenience, we have divided the proof of Theorem \ref{Res:MainCompactness} into a sequence of lemmas that in our opinion have their own inherent interest.

\begin{Lemme} If $v>0$ is fixed, then there exists $\tilde{M}:=M\bigcup_{i=1}^N M_{\infty,i}$ be a disjoint union of finitely many limit manifolds $(M_{\infty,i}, g_{\infty,i})=\lim_j(M,p_{i,j}, g)$, such that $I_{\tilde{M}}(v)=I_{M}(v)$. Moreover $I_{\tilde{M}}(v)$ is achieved by a generalized isoperimetric region.
\end{Lemme}
\begin{Dem}
It is a trivial to check that $I_{M}\geq I_{\tilde{M}}$. Observe that when $d_M(p_{ij}, p_{lj})\leq K$ for some constant $K>0$, we have $(M_{i,\infty}, g_{i,\infty})=(M_{l,j}, g_{l,\infty})$.  Thus we can restrict ourselves to the case of sequences diverging in different directions, i.e., $d_M(p_{i,j}, p_{l,j})\rightarrow+\infty$ for every $i\neq j$.  Then mimicking the proof of the preceding lemma, one can obtain $I_{M}\leq I_{\tilde{M}}$.
\end{Dem}

Here we give the proof of Theorem \ref{Res:MainCompactness}, which is an intrinsic generalization of \cite{Mor94} adapted to the context of manifolds with bounded geometry.\\
\begin{Lemme}\label{Lemma: SequenceStructure}[Bounded Volume and Area Sequence's Structure] Let $M^n$ be a complete Riemannian manifold with bounded geometry. Let $(\Omega_j)$ be a sequence of finite perimeter sets with $\P(\Omega_j)\leq A$ and $V(\Omega_j)\leq v$. Then there exist a subsequence $(\Omega_{j_k})$ that we rename by $(\Omega_k)$, a set $S\subseteq\N$, with $k\in S$, a natural number
$l\in\left(\N\setminus\{0\}\right)\mathring{\cup}\{+\infty\}$, $l$ sequences of points $(p_{ik})\in M^n$, $l$ sequences of radii $R_{ik}\rightarrow +\infty$, $i\in\{1, ..., l\}$, a sequence of volumes $v_i\in ]0,+\infty[$, and areas $A_i\in [0,+\infty[$ such that 
\begin{enumerate}[(I):]
           \item $v_i=\lim_{k\rightarrow+\infty}V(\Omega_{ik})$, where $\Omega_{ik}:=\Omega_k\cap B(p_{ik}, R_{ik})$,
           \item $0<\sum_{i=1}^{\infty}v_i=\bar{v}\leq v$,
           \item If $V(\Omega_k)\rightarrow v$, then $\bar{v}=v$,
           \item Set $A_i=\lim_{k\rightarrow+\infty}\P(\Omega_{ik})$, we get $$\lim_{k\rightarrow+\infty}\P(\mathring{\cup}_{i=1}^{\infty}\Omega_{ik})=\sum_{i=1}^{\infty} A_i=\bar{A}\leq A.$$
\end{enumerate}
\end{Lemme}
\begin{Rem} In (IV): $\bar{A}$, could be strictly less than $A$ even if $\P(\Omega_k)\rightarrow A$. However there are cases in which $\bar{A}$ is equal to $A$, as showed in Corollary \ref{Cor:amsequencestructure}.
\end{Rem}

\begin{Dem}[of Lemma \ref{Lemma: SequenceStructure}] This proof consist in finding a subsequence $\Omega_k$ of an arbitrary sequence $\Omega_j$ of bounded volume and area that could be decomposed as an at most countable union $\Omega_k=\bigcup _i^N \Omega_{ik}$, with $N$ possibly equal to $+\infty$. To this aim, we consider the concentration functions $$Q_{1,j}(R):=Sup_{p\in M}\{V(\Omega_j\cap B_M(p,R))\}.$$ Since Lemma 2.5 of \cite{NarAsian}    
 prevents evanescence in bounded geometry, the concentration-compactness argument of Lemma I.1 of \cite{Lions} or Lemma 2.5 of \cite{NarAsian}, provides a concentration volume $0<v_1\in[m_0(n,r=1, k,v_0, v),v]$.

Suppose the concentration of volumes occurs at points $p_{1j}$, a  ...
then there exists $E\subseteq\R$, with $|\R\setminus E|=0$, $S\subseteq\N$ such that for every $\epsilon_k\rightarrow 0$ there exists $R_{\epsilon_k}$ such that for every $R\in E\cap[R_{\epsilon_k}, +\infty[$ there exists $j_{k, R}$ satisfying the property that if $j\geq j_{k, R}$ then 
\begin{equation}
          |V(B(p_{k,R}, R)\cap\Omega_{j})-v_1|\leq\varepsilon_k. 
\end{equation} 
By the fact that the preceeding equation is true for every $R\in E\cap[R_{\epsilon_k}, +\infty[$, we can associate to every $k$ a radius $R'_{1k}$ satisfying
\begin{equation}
R'_{1,k+1}\geq R'_{1,k}+k,
\end{equation}
then by the coarea formula we get radii $R_{1k}\in[R'_{1,k}, R'_{1,k+1}]$ such that
\begin{equation}
\P(B(p_{k,R_{1k}}, R_{1k})\cap\Omega_{j})\leq A+\frac{v}{j}.
\end{equation}
We can set now $\Omega_{1k}:=B(p_{k,R_{1k}}, R_{1k})\cap\Omega_k$ for $k$ sufficiently large inside $S$. If $v_1=v$ then the theorem is proved with $l=1$. If $v_1<v_2$, we apply the same procedure to the domains $\Omega'_{1k}:=\Omega_k\setminus B(p_{k,R_{1k}}, R_{1k})$, observing that $V(\Omega'_{1k})\rightarrow v-v_1$, we obtain $R_{2k}$, $p_{2k}$,
$v_2$, $A_2$. If $v_1+v_2=v$, then we finish the construction and the theorem is proved in a finite number of steps. If this does not happen then we obtain at the $i$-th step a decomposition $\Omega_k=\Omega'_{ik}\mathring{\cup}_{l=1}^i\Omega_{ik}$. As $k$ goes to infinity we have
$V(\Omega_{ik})\rightarrow v_i$ and $V(\Omega'_{ik})\rightarrow v-\sum_{l=1}^iv_l$, if $V(\Omega_k)\rightarrow v$, and $\P(\partial\Omega'_{ik})\rightarrow A-\sum_{l=1}^iA_l$ where, $$A_l:=\lim_{k\rightarrow +\infty}\P((\partial\Omega_k)\cap B(p_{lk}, R_{lk})),$$ if $\P(\Omega_k)\rightarrow A$. In any case $\P(\Omega'_{ik})\leq A+i\frac{v}{k}$, so for $k$ big enough $\P(\Omega'_{ik})\leq 2A$. At this stage we found a sequence $(v_i)$ of positive volumes giving rise to a convergent series $\sum_{i=1}^{\infty}v_i=:\bar{v}\leq v$. To prove that if $V(\Omega_j)\rightarrow v$ then $\bar{v}=v$ we use the non evanescence Lemma 2.5 \cite{NarAsian}. Indeed $v_i\rightarrow 0$ because $\sum_{i=1}^{\infty}v_i$ is a convergent series, on the other end an application of  Lemma 2.5 of \cite{NarAsian} to $\Omega'_{ik}$ yields to
\begin{equation}
V(\Omega'_{ik}\cap B(p_{i,k+1}, R_{i,k+1}))=:v_{i+1,k}\geq\frac{V(\Omega'_{ik})^n}{\P(\Omega'_{ik})^n+1},
\end{equation}
but $v_{i+1}\geq v_{i+1,k}$ so we obtain $$\lim_{i\rightarrow+\infty}v_i=\lim_{i\rightarrow+\infty}v_{i+1}\geq\lim_{i\rightarrow+\infty}\limsup_{k\rightarrow +\infty}v_{i+1,k}\geq C(n,k,v_0,v)\frac{(v-\bar{v})^n}{A^n+1}>0.$$
This last inequality contradicts the fact that $v_i\rightarrow 0$, if $\bar{v}<v$. It remains to prove (V), but this is a simple consequence of the definitions.
\end{Dem}
\begin{Cor}\label{Cor:amsequencestructure}
Let $M^n$ be a complete Riemannian manifold with bounded geometry. If $(\Omega_j)$ is an almost minimizing sequence of volume $v$  then there exist a finite number $N=N(n, v_0, k, v)\in\N\setminus\{0\}$, a subsequence $(\Omega_{j_k})$ that we rename by $(\Omega_k)$, $k\in S\subseteq\N$, $N$ sequences of points $(p_{ik})_{i\in \{1,...,N\}}$, $N$ sequences of radii $(R_{ik})_{i\in \{1,...,N\}}\rightarrow +\infty$, when $k\rightarrow +\infty$, $N$ volumes $v_i\in ]0,+\infty[$, with $i\in\{1,..., N\}$, $N$ areas $A_i\in [0,+\infty[$, such that
\begin{enumerate}[(I):]
           \item $v_i=\lim_{k\rightarrow+\infty}V(\Omega_{ik})$, where $\Omega_{ik}:=\Omega_k\cap B(p_{ik}, R_{ik})$,
           \item $0<\sum_{i=1}^{N}v_i=\bar{v}\leq v$,
           \item If $V(\Omega_k)\rightarrow v$, then $\bar{v}=v$,
           \item $A_i=\lim_{k\rightarrow+\infty}\P(\Omega_{ik})$, moreover $\sum_{i=1}^{N} A_i=\bar{A}\leq A$
           \item If $\P(\Omega_k)\rightarrow A$, then $\bar{A}=A$.
           \item $N\leq\left[\frac{v}{v^*}\right]+1$, where $v^*=v^*(n, k, v_0)$ can be taken equal to $\bar{v}$ of Lemma 3.2 of \cite{EB}.
           In particular if $v<v^*$ then $N=1$.
\end{enumerate} 
\end{Cor}
\begin{Dem} The proof goes along the same lines of Theorem $1$ of \cite{NarAsian}. We prove just that $\bar{A}=A$. In view of Theorem 1 of \cite{FloresNardulli} the isoperimetric profile in bounded geometry is continuous so it is enough to prove (V) in the case of a minimizing sequence in volume $v$. Suppose, now, that $V(\Omega_k)=v$ and $A=I_M(v)$. From Lemma \ref{Lemma: SequenceStructure} we get a decomposition of $\Omega_k=\Omega'_k\mathring{\cup}\tilde{\Omega}_k$, where $\tilde{\Omega}_k:=\mathring{\cup}_{i=1}^{\infty}\Omega_{ik}$ (this infinite union in the case of a minimizing sequence is in fact a finite union) and $\Omega'_k:=\Omega_k\setminus\tilde{\Omega}_k$. It is easy to see that $V(\Omega'_k)=:v'_k\rightarrow 0$ and $V(\tilde{\Omega}_k)=:\tilde{v}_k\rightarrow 0$ as $k\rightarrow 0$, which implies
$$ I_M(v)=A=\lim_{k\rightarrow+\infty}I_M(\tilde{v}_k)\leq\lim_{k\rightarrow+\infty}\P(\tilde{\Omega}_k)=\bar{A}.$$
The second equality is due to the continuity of the isoperimetric profile, the remaining are easy consequences of the definitions. To obtain an estimate on $N$, one way to proceed is to show that there is no dichotomy for volumes less than some fixed $v^*=v^*(n, k, v_0)>0$. Assuming the existence of such a $v^*$, we observe that the algorithm produces $v_1\geq v_2\geq\cdots\geq v_N$.  Furthermore, $v_N\leq v^*\leq v_{N-1}$ because it is the first time that dichotomy cannot appear, which yields
\begin{equation}
v^*(N-1)\leq v_{N-1}(N-1)\leq\sum_{i=1}^{N-1} v_i\leq \sum_{i=1}^N v_i=v.
\end{equation}
Consequently
\begin{equation}\label{eq:Mainestimates}
N\leq\left[\frac{v}{v^*}+1\right]=\left[\frac{v}{v^*}\right]+1,
\end{equation}
where $v^*=v^*(n, k, v_0)$ can be taken equal to $\bar{v}$ of Lemma 3.2 of \cite{EB}.
On the other hand, one can construct examples such that for every $v$, there are exactly $N=\left[\frac{v}{v^*}\right]+1$ pieces. So in this sense, the estimate (\ref{eq:Mainestimates}) is sharp and it concludes the proof of the theorem.
\end{Dem}

Now we prove Theorem \ref{Res:MainCompactness}.

\begin{Dem}
By Gromov's compactness theorem and a cumbersome diagonalization process there is a subsequence converging in the multipointed Gromov-Hausdorff topology, by the results of \cite{Colding97Ann} volumes converges. Just assuming Gromov-Hausdorff convergence we don't know if perimeter converges, but the hypothesis of $C^0$-bounded geometry ensures the convergence of perimeters.
\end{Dem} 
\begin{Thm}\label{Thm:Continuity} If $M^n$ have $C^0$ asymptotically bounded geometry, then the isoperimetric profile $I_M$ is continuous on $[0, V(M)[$.
\end{Thm}
\begin{Dem} Take a sequence of isoperimetric regions $\Omega_i\subset M$, such that $V(\Omega_i)=v_i\rightarrow v$. Then apply generalized compactness to extract a subsequence $\Omega_k$ converging in the multipointed $C^0$ topology to a generalized isoperimetric region $\Omega$ of volume $v$. From the existence of the generalized isoperimetric region $\Omega$ we deduce that $I_M$ is upper semicontinuous as in Theorem 2.1 of \cite{NarFlores} and from lower semicontinuity of the perimeter with respect to the multipointed $C^0$ topology we get lower semicontinuity of $I_{\tilde{M}}=I_M$. 
\end{Dem}
\begin{Rem} Using Theorem \ref{Thm:Continuity}, we can prove easily Corollary $1$ of \cite{NarAsian} and Corollary $1$ of  \cite{NarFlores} without using Theorem $1$ of \cite{NarFlores}. Of course Theorem $1$ of \cite{NarFlores} is stronger than Theorem \ref{Thm:Continuity} and have its own interest independently of the use that we did in proving Corollary $1$ of  \cite{NarFlores}.
\end{Rem}
\section{Appendix: weaker assumptions on the main theorem}
  \begin{Def}[Flat distance modulo $\nu$ of Ambrosio-Kircheim] 
$${\dFZ}_{\nu}(T_1,T_2):=\inf\left\{\mathbf{M}(R)+\mathbf{M}(S):T_1-T_2=R+\partial S+\nu Q, R\in\mathbf{I}_k, S\in\mathbf{I}_{k+1}\right\}$$
\end{Def}
\begin{Def}[Intrinsic flat distance Lang-Wenger, Sormani-Wenger]
Let $X, Y$ be two metric spaces carrying a current structure. Let us define $${\dF}^{\nu}(X,Y):=\inf\left\{\dFZ(\varphi_1(X),\varphi_2(Y))\right\},$$
where the infimum is taken over all isometric embeddings $\varphi_1:X\rightarrow Z$ and $\varphi_2:Y\rightarrow Z$, and ${\dFZ}_{\nu}$ is the flat distance modulo $\nu$ of Ambrosio-Kircheim in $Z$.
\end{Def}
\begin{Def} [Pointed Intrinsic flat convergence Lang-Wenger, Sormani-Wenger]
Let $(X_n,d_n,p_n), (X,d,p)$ be two pointed metric spaces we say that they converges in the \textbf{pointed intrinsic flat topology}, if for every $R>0$ we have 
$$\dF(B_{X_n}(p_n, R),B_X(p,R))\rightarrow 0.$$ 
\end{Def}
\begin{Def}
Let $(X, d)$ be a metric space, given $l$ distinct points $p_1, ..., p_l\in X$, we call the $(l+2)$-tuple $(X, d, p_1, ..., p_l)$ an \textbf{$l$-pointed metric space}. In the sequel we will allow to take even a countable set of points of $X$ and we will speak of a \textbf{multipointed metric space} denoting by $(X, d, (p_j)_{j\in\N})$ or $(X, d, (p_j))$ when no confusion could arise. 
\end{Def}
\begin{Def}[Multipointed intrinsic flat convergence]
We say that a sequence of multipointed current metric spaces $(X_i, d_i, p_{i,1}, ..., p_{i,l},...)$ converges to the multipointed current metric space $(X_{\infty}, d_{\infty}, p_{\infty,1}, ..., p_{\infty,l},...)$, in the \textbf{multipointed intrinsic flat topology}, if for every index $j$, we have $$(X_i, d_i, p_{i,j})\rightarrow (X_{\infty}, d_{\infty}, p_{\infty,j}),$$ in the pointed intrinsic flat topology. 
\end{Def} 
\begin{Def} Let $(Z, d_Z)$ be a metric space, the \textbf{Hausdorff distance} between $A,B\subseteq Z$ is defined by 
$$d_{Haus}(A,B):=\left\{\varepsilon:A\subseteq\mathcal{U}_{\varepsilon}(B)\text{and} B\subseteq\mathcal{U}_{\varepsilon}(A)\right\},$$ 
where $\mathcal{U}_{\varepsilon}(A)$ denote the \textbf{$\varepsilon$-neighbourhood} of a subset $A\subseteq Z$, i.e., $\mathcal{U}_{\varepsilon}(A):=\{x\in Z: d_Z(x,A)\leq\varepsilon\}$.
\end{Def}
\begin{Def}[Gromov-Hausdorff distance]
Let $X, Y$ be two metric spaces. Let us define $$d_{\mathcal{GH}}(X,Y):=\inf\left\{d_{\mathcal{H}}^Z(\varphi_1(X),\varphi_2(Y))\right\},$$
where the infimum is taken over all isometric embeddings $\varphi_1:X\rightarrow Z$ and $\varphi_2:Y\rightarrow Z$, and $d_{\mathcal{H}}^Z$ is the Gromov-Hausdorff distance in $Z$.
\end{Def}
\begin{Def}[Gromov-Hausdorff convergence] A sequence of metric spaces $(X_n, d_n)$ is said to converge to a metric space $(X, d)$, if $d_{\mathcal{GH}}(X_n, X)\rightarrow 0$, when $n\rightarrow+\infty$. 
\end{Def}
\begin{Def}[Pointed Gromov-Hausdorff convergence] Let $(X_n, d_n, p_n)$ be a sequence of pointed metric spaces with base point $x_n\in X_n$.  
The sequence $(X_n, d_n, p_n)$ is said to converge to a pointed metric space $(X,d,x)$, if for every $R>0$ the sequence of closed balls $\bar{B}(x_n, R)$, with the induced metrics converges to the closed ball $\bar{B}(x,R)\subseteq X$ in the Gromov-Hausdorff topology. 
\end{Def}
\begin{Rem} It is easy to see, from the definitions, that in general, if two metric spaces have $0$ Gromov-Hausdorff distance they do not be necessarily isometric. For example is zero the G-H distance between a metric space and any dense subspace of it, hence G-H limits are not unique in general. However, if $X$, $Y$ are two compact metric spaces then they are isometric iff $d_{\mathcal{GH}}(X,Y)=0$. This ensures uniqueness of limits of sequences of compact spaces that converges to another compact space. Observe here that a sequence of compact metric spaces could have noncompact G-H limits.  
\end{Rem}
\begin{Def}
A metric space $(X, d)$ is called \textbf{proper}, if for each point $x_0\in X$, the distance function to $x_0$, $x\mapsto d(x_0,x)$ is a proper function $d(x_0, \cdot):X\rightarrow\R$, i.e., if for each point $x_0\in X$ and $R>0$ the \textbf{closed ball} $\bar{B}(x_0, R):=\{x\in X:d(x_0,x)\leq R\}$ is compact.  
\end{Def}
\begin{Def}[Multipointed intrinsic Gromov-Hausdorff convergence]
We say that a sequence of multipointed metric spaces $(X_i, d_i, p_{i,1}, ..., p_{i,l},...)$ converges to the multipointed metric space $(X_{\infty}, d_{\infty}, p_{\infty,1}, ..., p_{\infty,l},...)$, in the \textbf{multipointed Gromov-Hausdorff topology}, if for every $j$ we have $$(X_i, d_i, p_{i,j})\rightarrow (X_{\infty}, d_{\infty}, p_{\infty,j}),$$ in the pointed Gromov-Hausdorff topology. 
\end{Def}
\begin{Rem} We observe that if the sequence $X_n$ satisfies a uniform bounded geometry condition as in definition \ref{Def:BoundedGeometry}, then intrinsic flat convergence implies Gromov-Hausdorff convergence.
\end{Rem}
\begin{Thm} The perimeter is lower semicontinuous with respect to intrinsic flat convergence. The volume is continuous with respect to intrinsic flat convergence.
\end{Thm}
For the definitions needed in what follows we refer to \cite{LangWenger}.
\begin{Def}[\cite{LangWenger}]
We say that a sequence $(T_j)$ in $\textbf{I}_{loc,m}(Z)$ converges in the local flat topology to a current $T\in\textbf{I}_{loc,m}(Z)$ if for every bounded closed set $B\subseteq Z$ there is a sequence $(S_j)$ in $\textbf{I}_{loc,m+1}(Z)$ such that
$$\left(||T-T_j-\partial S_j||+||S_j||\right)(B)\rightarrow 0.$$
\end{Def}
\begin{Thm} Let $X$ be a complete metric space such that $[X]$ induces a current in $\textbf{I}_{loc,n}(X)$, let $T_k\in\textbf{I}_{loc,n}(X)$ a sequence. Suppose that $X$ is doubling or satisfy a curvature dimension condition with dimension $n$, satisfy a $1$-local Poincar\'e inequality, an Euclidean type isoperimetric inequality for small volumes, and 
\begin{equation}
\sup_n\{||T_k||(X) + ||\partial T_k||(X)]\}< C.
\end{equation} 
Then there exist a subsequence $(k_j)$, a complete multipointed metric space $Z=(\mathring{\cup}_iZ_i, z_i\in Z_i,d_Z)$, $\varphi_{ij} : X\hookrightarrow Z$ such that $\varphi_{ij}(x_{k_j})\rightarrow z_i$ and $(\varphi_{ij\#}(T_{k_j})\subseteq Z_i)$ converges in the local flat topology of $Z_i$ to some $T\in \textbf{I}_{loc,n}(Z_i)$.
\end{Thm}
\begin{Dem} 
\begin{enumerate}[(i):]
          \item The arguments of Lemma 2.2 applies mutatis mutandis and this produces a sequence of sequences of points $(p_{j})_i$ and volumes $v_i$
          \item Apply Theorem 1.1 of \cite{LangWenger} together with a cumbersome diagonal argument  to obtain limit spaces $Z_i$ and the theorem follows easily. 
\end{enumerate}        
\end{Dem}
\begin{Rem} The preceding Theorem is a generalization of Theorem \ref{Res:MainCompactness} in the more general setting of the theory of Sormani, Lang and Wenger.
\end{Rem}
\newpage 
      \markboth{References}{References}
      \bibliographystyle{alpha}
      \bibliography{these}
      \addcontentsline{toc}{section}{\numberline{}References}
      \emph{Stefano Nardulli\\ Departamento de Matem\'atica\\ Instituto de Matem\'atica\\ UFRJ-Universidade Federal do Rio de Janeiro, Brazil\\ email: nardulli@im.ufrj.br\\ \\ Abraham Enrique Mu\~noz Flores\\ Ph.D student\\Instituto de Matem\'atica\\UFRJ-Universidade Federal do Rio de Janeiro, Brazil\\email: abrahamemf@gmail.com} 
\end{document}